\documentclass[letterpaper,11pt]{article}
\usepackage{amsmath, amsthm, amssymb}    
\usepackage{bridges}			
\usepackage{graphicx}			
\usepackage[colorlinks=true, urlcolor=blue, citecolor=black, linkcolor=black]{hyperref}  
\usepackage{subcaption}			

\usepackage{amssymb}
\usepackage{graphicx} 
\usepackage{enumerate}
\usepackage{mathtools}
\usepackage{color,soul}
\usepackage{cite}
\usepackage{tikz}
\usetikzlibrary{matrix}
\usepackage{caption}
\usepackage{amsmath, amsthm}
 \usetikzlibrary{cd}
 \usetikzlibrary{decorations.pathreplacing}
 \usepackage{tikz,calc}
\usepackage{color}

\usepackage{multicol,todonotes}

\newcommand{\rpt}{\mathbb{R}\text{P}^2}

\newcommand{\pcross}{
\begin{tikzpicture}[baseline=-2.75, scale=.2]
\draw[-{Stealth[ length=1.25mm, width=1.25mm]},thick ](-1,-1)--(1,1);
\draw[-{Stealth[ length=1.25mm, width=1.25mm]},thick ](-.25,.25)--(-1,1);
\draw[thick ](.3,-.3)--(1,-1);
\end{tikzpicture}}

\newcommand{\ncross}{
\begin{tikzpicture}[baseline=-2.75,scale=.2]
\draw[-{Stealth[ length=1.25mm, width=1.25mm]},thick ](1,-1)--(-1,1);
\draw[-{Stealth[ length=1.25mm, width=1.25mm]},thick ](.25,.25)--(1,1);
\draw[thick ](-.3,-.3)--(-1,-1);
\end{tikzpicture}}

\newcommand{\swap}{
\begin{tikzpicture}[baseline=-2.75, scale=.2]
\draw[-{Stealth[ length=1.25mm, width=1.25mm]},thick ](-1,-1)--(1,1);
\draw[-{Stealth[ length=1.25mm, width=1.25mm]},thick ](1,-1)--(-1,1);
\draw(0,0) circle (12pt);
\end{tikzpicture}}

\newcommand{\tbar}{
\begin{tikzpicture}[baseline=-2.75, scale=.2]
\draw[-{Stealth[ length=1.25mm, width=1.25mm]},thick ](0,-1)--(0,1);
\draw[thick ](-.35,0)--(.35,0);
\end{tikzpicture}}

\newcommand{\halfcirc}{
\begin{tikzpicture}[baseline=-2.75, scale=.2]
\draw[fill] (0,0)-- (90:4ex) arc (90:270:4ex) -- cycle ;
  \draw (0,0) circle (4ex);
\end{tikzpicture}}

\newcommand{\fullcirc}{
\begin{tikzpicture}[baseline=-2.75, scale=.2]
  \fill (0,0) circle (4ex);
  \draw (0,0) circle (4ex);
\end{tikzpicture}}

\newcommand{\theoremname}{Theorem:}

\title{Danceability of Twisted Virtual Knots}
\author{Sol Addison\textsuperscript{1}, Nancy Scherich\textsuperscript{2}, and Lila Snodgrass\textsuperscript{3}
\vspace{10pt}\\
Dept. of Mathematics and Statistics, Elon University, Elon, North Carolina, USA\\
\textsuperscript{1} saddison3@elon.edu, 
\textsuperscript{2} nscherich@elon.edu, 
\textsuperscript{3} lsnodgrass@elon.edu} 

\date{}					
\begin{document}
\maketitle

\vspace{-.5mm}
\begin{abstract}
 Over the years, several Bridges papers have delved into the concept of \emph{danceability} of a knot diagram. Inspired by dancing on non-orientable surfaces, in this paper, we expand danceability to twisted virtual knot diagrams. This paper is accompanied by a Math-Dance video which can be found in \cite{VID}.\end{abstract}

\vspace{-3mm}
\section*{Introduction to Knot Diagrams and Danceability}

\begin{figure}
      \centering
      \begin{picture}(450,50)
      \put(3,8){\includegraphics[scale=.5]{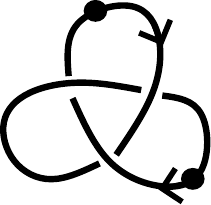}}
      \put(-2,0){(a)}
      \put(80,0){\includegraphics[scale=.4]{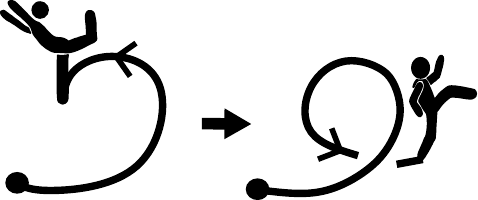}}
      \put(60,0){(b)}
      \put(180,-10){\includegraphics[scale=.19]{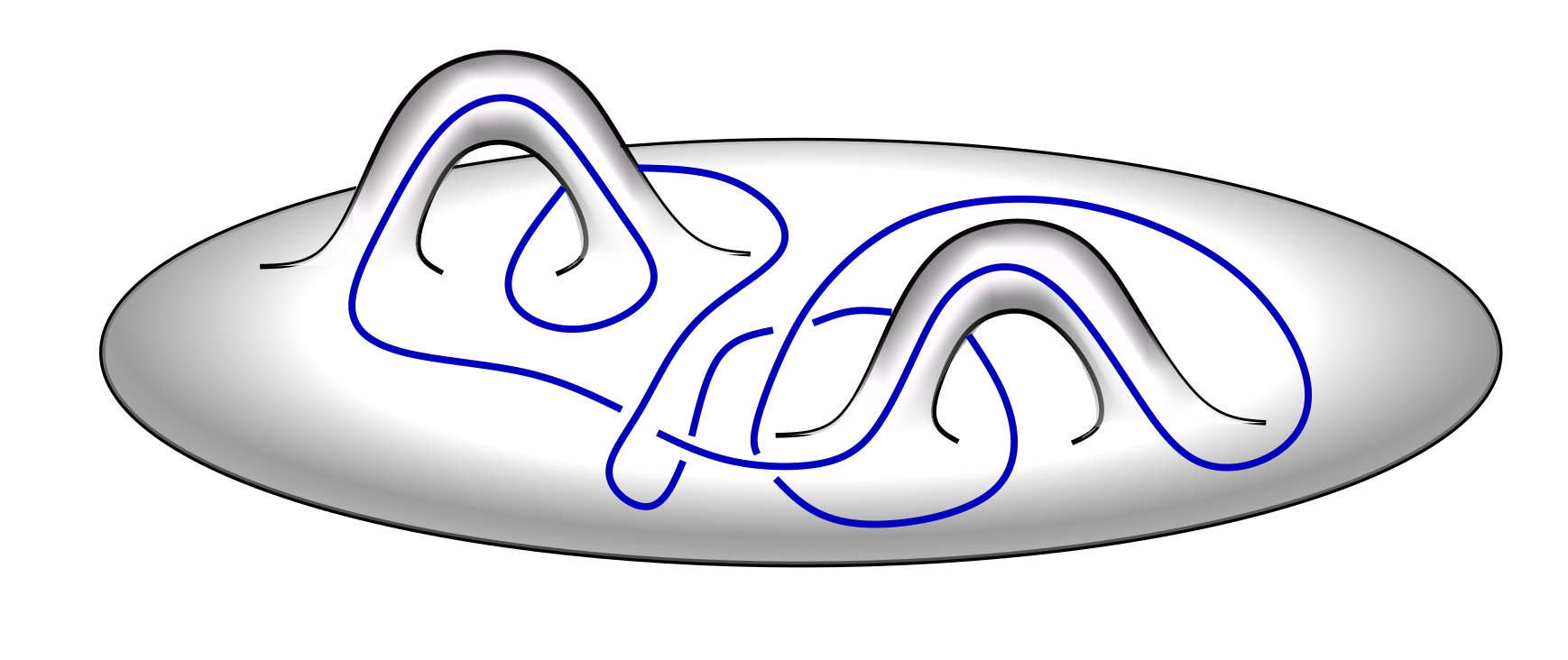}}
      \put(173,0){(c)}
    \put(350,0){\includegraphics[scale=.4]{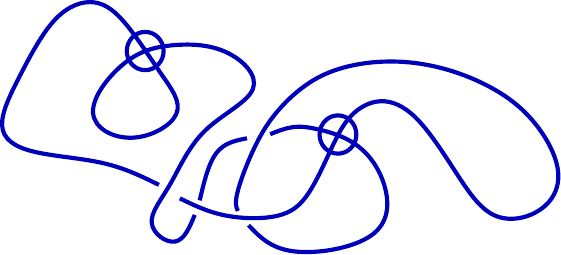}}
      \end{picture}
      \caption{(a) A 2-danceable trefoil. The two dots denote the starting positions of the two dancers. (b) A visualization of the over-first rule as a dancer's path traces out an arc. (c) A knot diagram on a genus 2 surface and its corresponding planar projection drawn as a virtual knot diagram.}
      \label{fig:over_first_rule}
  \end{figure}
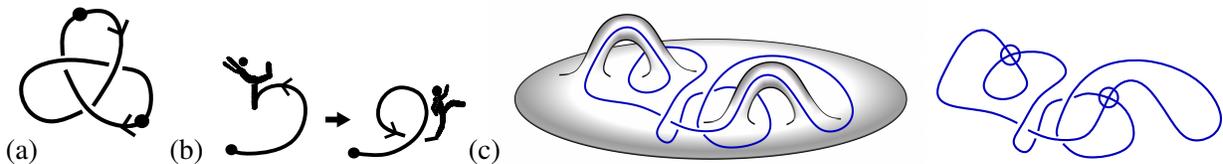

A \emph{knot} is a closed curve in $\mathbb{R}^3$ with no self intersections. From a topological perspective, knots are made of an extremely thin and flexible material; if you bend, stretch, or rearrange a knot without cutting it, a topologist would say it is the same knot. 
A \emph{knot diagram} is a type of planar projection of a knot with no triple points and where the over and under strands of a crossing are denoted in the diagram as $\pcross$ or $\ncross$. These crossings are known as \emph{classical} crossings. When a knot diagram consists solely of classical crossings, it is referred to as a classical knot. You can see an example of such a knot diagram in Figure \ref{fig:over_first_rule}(a).
Inspired by the knotted paths of dancers along a stage, one can study the \emph{danceability} of a knot diagram.

 \vspace{.2cm}

\noindent \textbf{Definition 1.} [Schaffer \cite{KS}, A.-S.-S. \cite{ASS}]
For $n\geq 2$, an oriented knot diagram is \emph{n-danceable} if there exists $n$ initial points on the knot diagram so that one dancer starts at each of the $n$ initial points, for a total of $n$ dancers.
  Each dancer travels in the direction of the pre-chosen orientation of the diagram, and stops dancing when they reach the next initial point. 
   The speed that each dancer travels can vary and be chosen so that the simultaneous tracing of the diagram by all $n$ dancers follows the over-first rule at every (classical) crossing.

\noindent \textbf{Over-first rule:} At every crossing, a dancer must pass through the crossing on the over-strand first before a dancer (possibly the same dancer) can cross on the under-strand. If a dancer reaches a crossing on an under-strand, they must wait until another dancer passes on the over-strand before continuing forward. See Figure \ref{fig:over_first_rule}(b). 

     \vspace{.2cm}
    

  \vspace{.2cm}

 Karl Schaffer first introduced the concept of danceability  in his foundational 2021 Bridges article \cite{KS}.  
 In another 2024 Bridges article \cite{ASS}, A.-S.-S.  extended Schaffer's ideas by formalizing the definition of an $n$-danceable knot, Definition 1, and studied the induced danceability knot invariant by taking the minimum number of dancers needed to dance any projection of a knot.
 A.-S.-S.-Sullivan proved that there is flexibility as to what crossing restriction rule one can use in the definition of $n$-danceable and still induce the same knot invariant \cite{ASSS}.
 As such,  this paper chooses to use the ``over-first" rule as the crossing restriction, while Schaffer's and A.-S.-S.'s original definitions used the ``under-first" rule, see \cite[Corollary 2.5]{ASSS}  for more details.

Interesting things happen if we draw knot diagrams on more complicated surfaces instead of just a plane. For example, Figure \ref{fig:over_first_rule}(c) depicts a knot diagram drawn on a genus two surface. Because of the topology of the surface, a strand can look like it crosses over another strand; however, this is a consequence of the shape of the surface and not actually a crossing of the diagram. This new type of crossing is called a virtual crossing and is denoted by \swap \cite{Kauf}. 
A \emph{virtual knot diagram } is a planar drawing of a knot that can have classical crossings $\pcross$ and $\ncross$, in addition to virtual crossings $\swap$. 
A.-S.-S.-Sullivan showed that the notion of $n$-danceable extends to virtual knot diagrams by adding a virtual crossing restriction rule that dancers must satisfy when dancing across a virtual crossing, see \cite[Definition 4.1]{ASSS}.
Dancers dancing a virtual knot diagram are still dancing on a plane (the stage is still flat), but the diagram they are dancing is a planar picture of a virtual knot which originally lays on a more complicated surface.

So far, the ideas of danceability have only been applied to knot diagrams that are drawn on \emph{orientable} surfaces (planes, spheres, genus $g$-surfaces, etc.). The goal of this paper is to explore how a dancer might dance  knot diagrams on  \emph{non-orientable} surfaces, called \emph{twisted virtual knot diagrams}.

\section*{What are Twisted Virtual Knots?}

\begin{figure}[ht]
    \begin{center}
    \hspace{5mm}
    \begin{picture}(130,40)
    \put(10,0){\includegraphics[scale=.11]{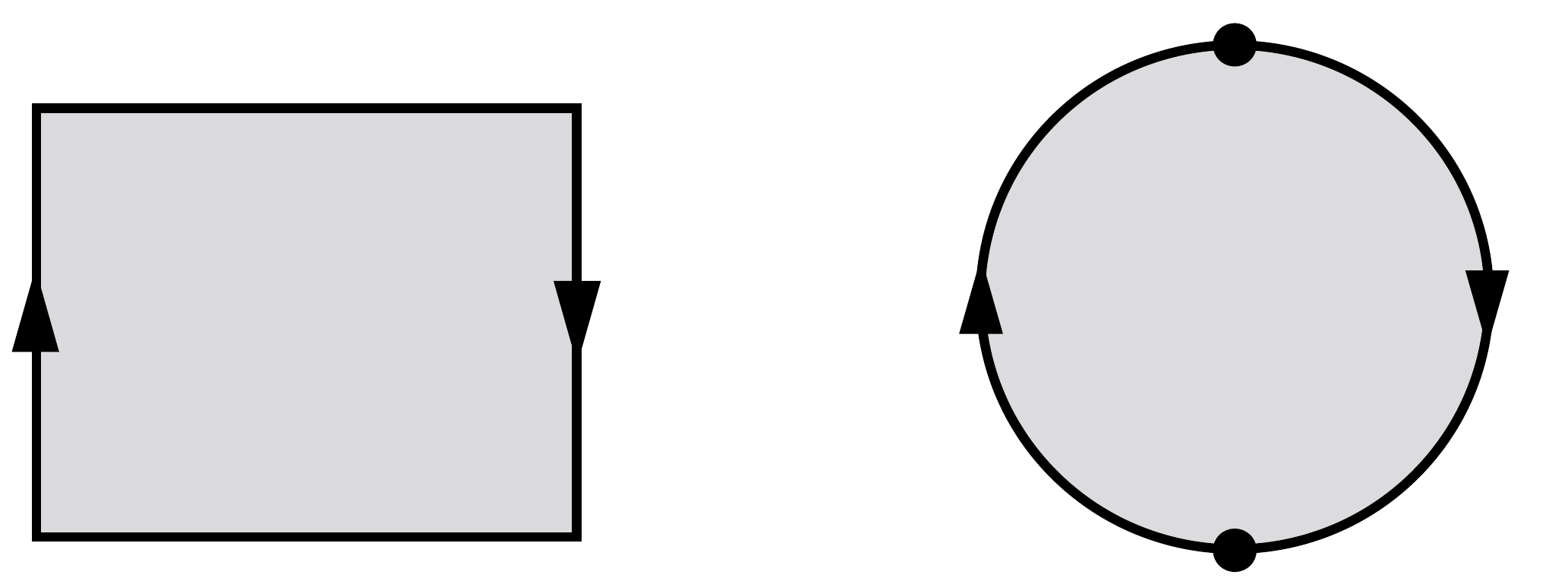}}
    \put(25,15){$\mathcal{M}$}
    \put(-5,0){(a)}
     \put(87,15){$\rpt$}
     \put(70,0){(b)}
    \end{picture}\begin{picture}(150,45)
   \put(0,0){
    \includegraphics[scale=.2]{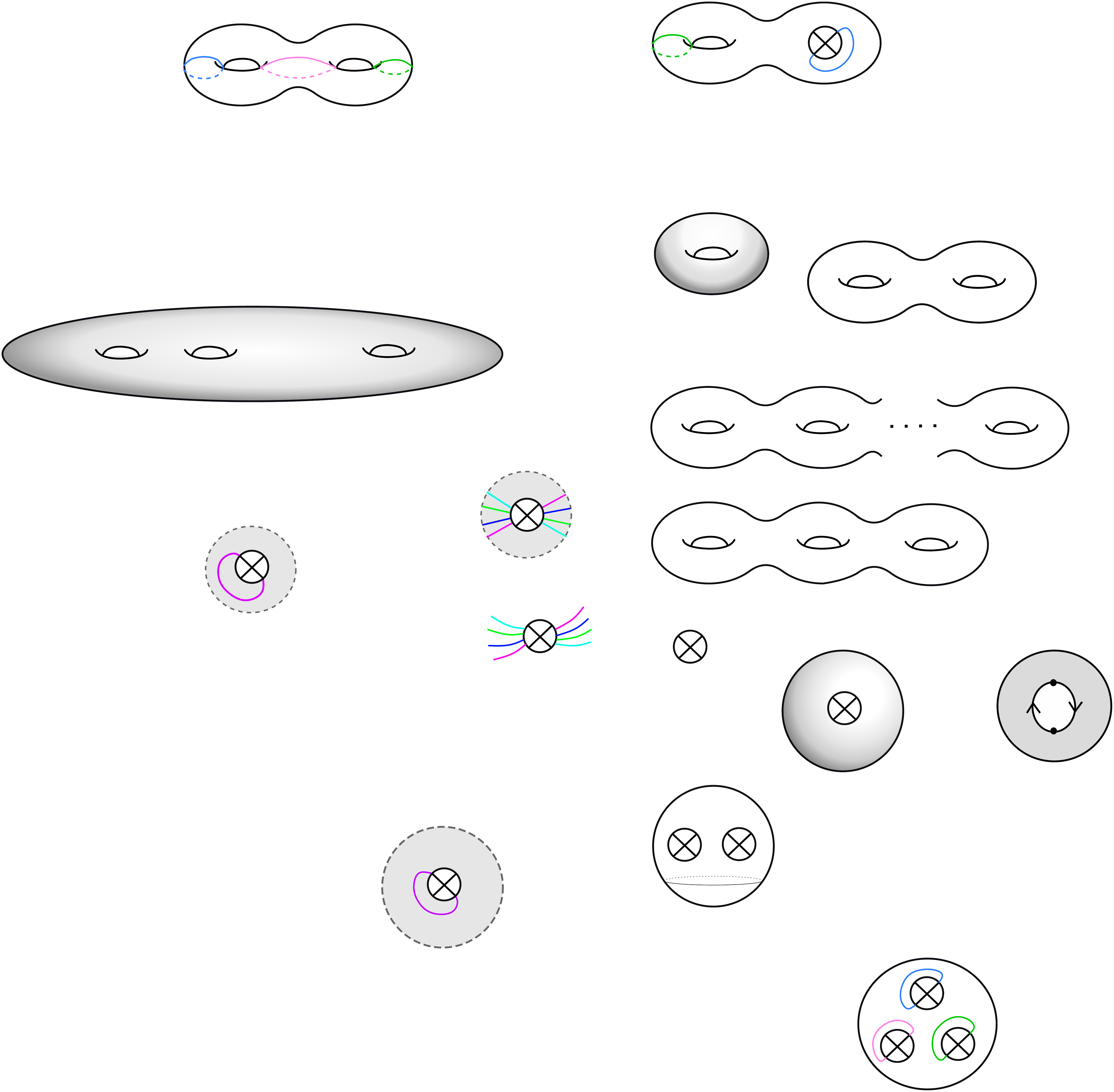}\hspace{5mm}}
    \put(75,0){\includegraphics[scale=.4]{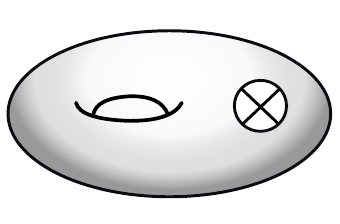}}
    \put(-7,0){(c)}
    \put(65,0){(d)}
    \end{picture}
    \hspace{2mm}
    \begin{picture}(150,45)
    \put(-10,0){(e)}
    \put(0,0){\includegraphics[width=0.27\linewidth]{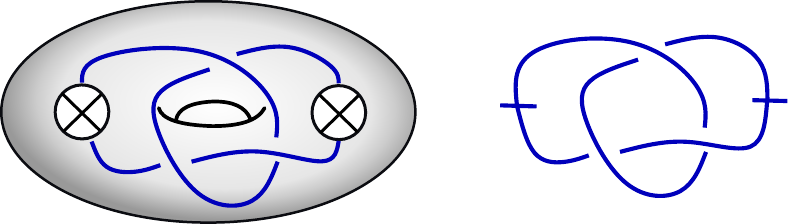}}
    \end{picture}
    \end{center}
    \caption{(a-b) Cut-and-paste diagrams for the Möbius band and  $\rpt$ (c) $\rpt$ minus a disc, also called a cross-cap.  (d) A torus connect sum  $\rpt$ drawn as a cross-cap on the torus. (e) A knot diagram on a torus with two cross caps and its corresponding planar projection drawn as a twisted virtual knot diagram. The two tick marks indicate twist bars.}
    \label{fig:exampleSurfaces}
    
\end{figure}
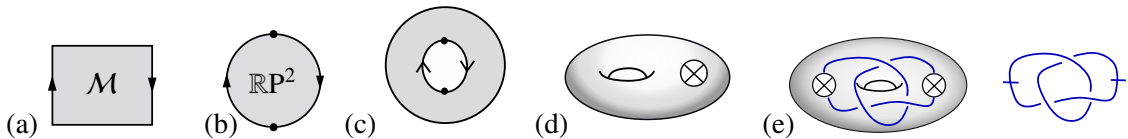

A  surface is \emph{non-orientable} if it is impossible to make a consistent choice of normal vector at every point on the surface (i.e. it is not possible to consistently tell the difference between clockwise and counterclockwise directions on the surface). 
The simplest example of an non-orientable surface is the Möbius band-- a rectangular strip of paper with opposite ends glued together with a  half twist, see Figure \ref{fig:exampleSurfaces}(a). 
The Real projective plane, denoted $\rpt$, is another example of a non-orientable surface and is shown in Figure \ref{fig:exampleSurfaces}(b).

One can glue together two surfaces using a \emph{connect sum}  operation, which consists of removing a disc from both surfaces and then gluing the two surfaces together along the circular boundary of the removed discs.
The classification of surfaces that states all (compact) non-orientable surfaces are a (compact) orientable surface connect sum a finite number of $\rpt$'s. 
To connect sum a surface $S$ with an $\rpt$, one must remove a disc from $\rpt$, as in Figure \ref{fig:exampleSurfaces}(c), and then glue that circular boundary to a circular boundary of $S-disc$.
Diagrammatically to denote a connect sum with an $\rpt$ we draw a $\otimes$ symbol\footnote{It is unfortunate that the established symbol $\otimes$ for a cross-cap looks very similar to the established symbol for a virtual crossing $\swap$. Cross-caps and virtual crossings are not related and one must use context clues to distinguish the symbol usage.}, called a cross-cap, on the surface where the $\rpt$ was connected. See Figure \ref{fig:exampleSurfaces} (d-e) for examples of cross-caps on surfaces.

Just as virtual knot theory is induced from knot diagrams drawn on oriented surfaces, \emph{twisted virtual knot theory} is induced from drawing knot diagrams on non-orientable surfaces \cite{Bou}.
Twisted virtual knots are diagrammatically generated by classical crossings $\pcross$ and $\ncross$, virtual crossings $\swap$, along with a new symbol $\tbar$ called a \emph{twist bar}. The twist bar indicates when a strand has passed through a cross-cap, see Figure \ref{fig:exampleSurfaces}(e) for an example.

\section*{Dancing Twisted Virtual Knots}


We can define different danceability numbers by imposing different rules for the dancers to adhere to 
at classical and virtual crossings (as in \cite[Definition 4.1]{ASSS}). In this section, we explore how a dancer can pass through a twist bar on a knot diagram. In the following definition for dancing twisted virtual knots, we have elected to put no restrictions on how a dancer may pass through a virtual crossing (this agrees with the unrestricted rule from \cite{ASSS}) 
in an effort to focus on the new twist bar restrictions.


A twist bar on a knot diagram represents when an arc of a knot has passed through a cross-cap. 
If a dancer were to walk through a cross-cap, they would instantaneously come out as the mirrored version of themselves. Since it is not physically possible to create a perfect mirror of a dancer, it is a choreographic question of how to represent such an orientation shift.
We choose to represent a ``mirrored" dancer as a dancer traveling along the oriented path backwards.
As such, dancers can be forward-facing or backward-facing.




\noindent \textbf{Definition 2.}
For $n,k>0$, an oriented twisted virtual knot diagram is \emph{forward $(n,k)$-danceable} (resp. \emph{matching  $(n,k)$-danceable}) if there are $n$ initial points on the diagram satisfying conditions 1-3 at every classical crossing as well as the respective fourth condition at a twist bar:
\vspace{-.2cm}
\begin{enumerate}
    \item One dancer starts at each of the $n$ initial points, for a total of $n$ dancers. (These initial points partition the diagram into $n$ oriented danceable segments called ``paths").
    \item Each dancer travels in the direction of the pre-chosen orientation of the diagram and
stops dancing when they have traversed $k$ paths.


    \item  All dancers start with a pre-chosen facing. At a twist bar, the dancer must flip their facing transitioning from forward-facing to backward-facing and vice versa.

    \item The speed that each dancer travels can vary and can be chosen so that the simultaneous tracing of the diagram by all $n$ dancers follows the respective rule:

\noindent \underline {Forward rule:} All dancers must start and end as a forward-facing dancer while adhering to the over-first rule at every classical crossing.

\noindent \underline {Matching rule:} Each initial point is designated a facing (either forward or backward) and a dancer starting at that initial point starts dancing in the designated facing. No dancer can end on an initial point whose designated facing does not agree with the current facing of the dancer. 






\end{enumerate}

Diagrammatically, we denote a forward-facing designation on an initial point with a $\fullcirc$ and the backward-facing initial point with a $\halfcirc$. 
Each dancer's current facing is denoted by a full circle (forward-facing dancer) or a half circle (backward-facing dancer) with different colors distinguishing each dancer. 
For the paths traversed by the dancers, a solid line denotes the path traversed by a forward-facing dancer, and a dotted line denotes a backward-facing dancer's path.

\begin{figure}
    \centering
    \begin{picture}(420,75)
    \put(20,0){\includegraphics[width=0.65\linewidth]{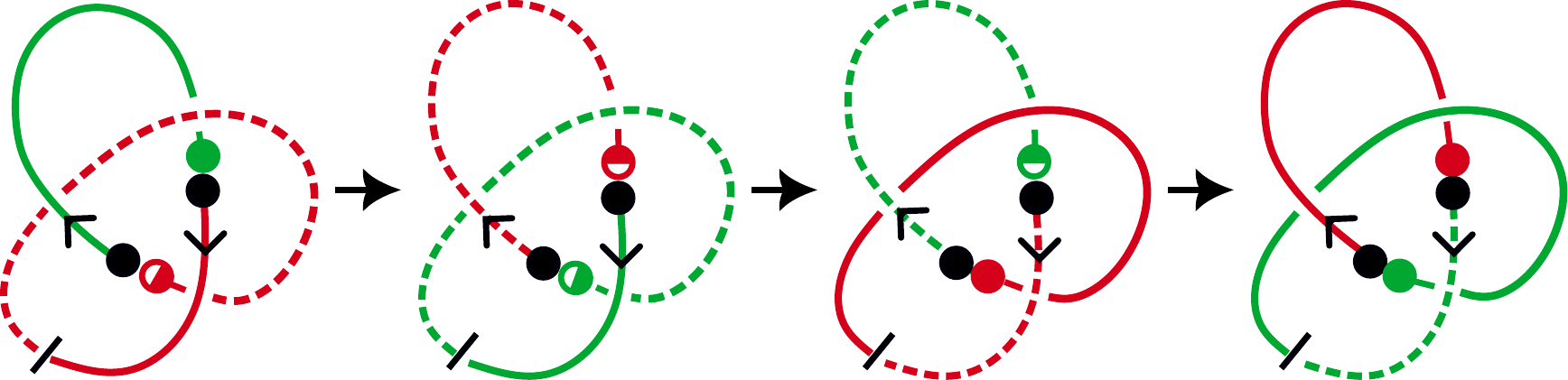}}
    \put(0,0){(a)}
    \put(360,0){\includegraphics[width=0.13\linewidth]{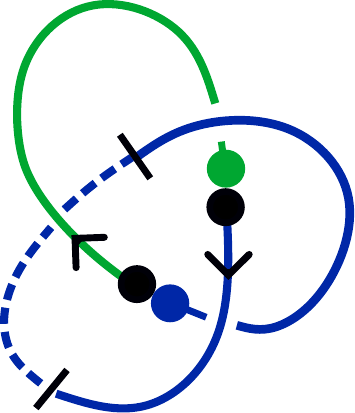}}
    \put(340,0){(b)}
    \end{picture}
    \caption{(a) The trefoil with one twist bar danced by 2 dancers following the forward rule. Each diagram is one iteration of the dancers walking a path. The dancers are shown at the end of their path. (b) The trefoil with two twist bars danced by 2 dancers following the forward rule with $k=1$.}
    \label{fig:trefoil_forward}
\end{figure}

 Figure \ref{fig:trefoil_forward}(a) shows an example of a trefoil knot with one twist bar danced with two dancers following the forward rule. 
The trefoil is divided into two distinct paths which the dancers must each traverse twice (for a total of $k = 4$ paths traveled) such that both dancers are facing forward at the end of their dance phrase. In essence, both dancers traversed the entire knot twice. 
This is an example of $k>n$ where all dancers traverse the entire knot and pass their initial starting point to continue dancing an additional $k-n$ paths. 

\begin{figure}
    \centering
    \begin{picture}(100,20)
    \put(0,0){\includegraphics[scale=.27]{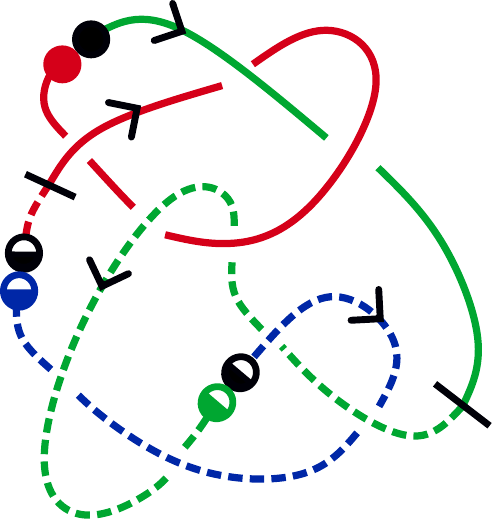}}
    \put(-5,-5){(a)}
     \end{picture}   \hspace{1cm} \begin{picture}(100,100)
    \put(0,0){\includegraphics[scale=.27]{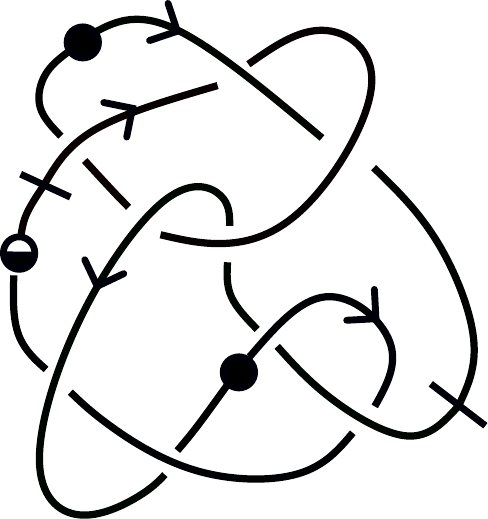}}
    \put(-5,-5){(b)}
     \end{picture}
    \caption{The same 9-crossing knot danced by 3 dancers following the matching rule with different initial facings. (a) Requires $k=1$. (b) Requires $k=3$. (iterations not shown)}
    \label{fig:9_crossing}
\end{figure}

In Figure \ref{fig:trefoil_forward}(b), we see a trefoil with two twist bars danced with two dancers following the forward rule. Since both twist bars are on the same path, the dancer switches facings twice, so the knot can be danced with $k=1$.
The two diagrams in Figure \ref{fig:9_crossing} are the same 9-crossing knot with two twist bars and with initial positions in the same locations. The facings in the diagrams agree for two initial positions and disagree at the third. Following the matching rule, the diagram in (a) can be danced with three dancers with $k=1$, whereas the diagram in (b) can be danced with three dancers but requires $k=3$.


 \vspace{-4mm}

\section*{Future Directions} The  motivation for Definition 2 focused on choreographic implications of dancing through a twist-bar. However, there are interesting mathematical consequences worth exploring. One can define and study the knot invariant induced by Definition 2 by taking the minimum $n$
 needed to dance a knot over all choices of  diagrams representing the knot.  
Another future direction is to study a  crossing rule  inspired by the idea of ``inversions" introduced in the  movie \textit{Tenet} directed by Christopher Nolan. In \textit{Tenet}, the ability to reverse entropy causes time inversions. For everyone else, time is still proceeding in a forward direction, but from the inverted person's perspective, the rest of the world is moving in reverse.  By the retrograde trick introduced in \cite{ASSS}, time-reversing a dance following the over-first rule gives an under-first dance. 
We conceived the \textit{Tenet rule} that requires a forward-facing dancer to adhere to the over-first rule, while a backward-facing dancer must adhere to the under-first rule. The Tenet rule creates new complications at mixed crossings between one forward-facing and one backward-facing dancer. This rule requires more clarity before it can be  implemented.

 \vspace{-4mm}

\section*{Acknowledgments}
This project is part of an undergraduate student research experience at Elon University for the first and third author, led by the second author.
This project is funded in part by the Lumen Prize.

{\setlength{\baselineskip}{13pt} 
\raggedright				

\end{document}